\documentclass{amsart}
\usepackage{amsmath,amssymb,amsthm}
\usepackage[longnamesfirst]{natbib}
\usepackage{setspace}
\usepackage{enumerate}         

\newtheorem{thm}{Theorem}[section]

\newtheorem{lem}[thm]{Lemma}
\newtheorem{prop}[thm]{Proposition}
\theoremstyle{definition}
\newtheorem{defn}[thm]{Definition}

\newtheorem{ass}[thm]{Assumption}
\theoremstyle{remark}
\newtheorem{rem}[thm]{Remark}

\numberwithin{equation}{section}

\newcommand{\Pt}{(P_t)_{t\geq0}}
\newcommand{\Xt}{(X_t)_{t\geq0}}

\newcommand{\norm}[1]{\left\Vert#1\right\Vert}

\newcommand{\set}[1]{\left\{#1\right\}}

\newcommand{\RR}{\mathbb{R}}

\newcommand{\PP}{\mathbb{P}}
\newcommand{\CC}{\mathbb{C}}
\newcommand{\NN}{\mathbb{N}}
\newcommand{\EE}{\mathbb{E}}
\newcommand{\cU}{\mathcal{U}}

\newcommand{\cF}{\mathcal{F}}

\newcommand{\cQ}{\mathcal{Q}}
\newcommand{\cR}{\mathcal{R}}

\newcommand{\cL}{\mathcal{L}}
\newcommand{\cT}{\mathcal{T}}
\newcommand{\Rplus}{\mathbb{R}_{\geqslant 0}}

\newcommand{\pd}[2]{\frac{\partial #1}{\partial #2}}

\newcommand{\scal}[2]{\left\langle{#1},{#2}\right\rangle}

\renewcommand{\Re}{\mathrm{Re}}
\renewcommand{\Im}{\mathrm{Im}}

\newcommand{\wt}[1]{{\widetilde{#1}}}
\newcommand{\Ex}[2]{\mathbb{E}^{#1}\left[#2\right]}                     
\newcommand{\Excond}[3]{\mathbb{E}^{#1}\left[\left.#2\right|#3\right]}  
\newcommand{\wh}[1]{{\widehat{#1}}}
\newcommand{\id}{\operatorname{id}}

\title{Affine processes are regular}
\author{Martin Keller-Ressel}
\address{ETH Z\"urich, D-Math, R\"amistrasse 101, CH-8092 Z\"urich, Switzerland}
\email{kemartin@math.ethz.ch}

\author{Walter Schachermayer}
\address{University of Vienna, Faculty of Mathematics, Nordbergstrasse 15, 1090
Vienna, Austria}
\email{walter.schachermayer@univie.ac.at}

\author{Josef Teichmann}
\address{ETH Z\"urich, D-Math, R\"amistrasse 101, CH-8092 Z\"urich, Switzerland}
\email{jteichma@math.ethz.ch}

\thanks{The first and the third author gratefully acknowledge the support from
the Austrian Science Fund (FWF) under grant Y328 (START prize).}
\thanks{The second and third author gratefully acknowledge the support from the
Austrian Science Fund (FWF) under grant P19456, from the Vienna Science and
Technology Fund (WWTF) under grant MA13 and by the Christian
Doppler Research Association (CDG)}

\begin{document}

\begin{abstract}
We show that stochastically continuous, time-homogeneous affine
processes on the canonical state space $\Rplus^m \times \RR^n$ are
always regular. In the paper of \citet{Duffie2003} regularity was
used as a crucial basic assumption. It was left open whether this regularity
condition is automatically satisfied, for stochastically continuous affine
processes. We now show that the regularity assumption is
indeed superfluous, since regularity follows from stochastic continuity and the
exponentially affine behavior of the characteristic function. For the proof we
combine classic results on the differentiability of transformation semigroups
with the method of the moving frame which has been recently found to be useful
in the theory of SPDEs.
\end{abstract}

\keywords{affine processes, regularity, characteristic function, semiflow. MSC
2000: 60J25, 39B32}

\maketitle

\section{Introduction}

A Markov process $ X $ taking values in $D:=\RR_{\geq 0}^m \times
\RR^n $ is called \emph{affine} if there exist $ \mathbb{C}$-valued
functions\footnote{We are using upper case $\Phi$ but lower case $\psi$ for
consistency
with the notation of \citet{Duffie2003}; see Remark~\ref{rem:phi}
for a detailed discussion.} $\Phi(t,u)$ and $\psi(t,u)$ such that
\begin{equation}\label{Eq:affine_prop1}
\EE^x\left[e^{\scal{X_t}{u}}\right] = \Phi(t,u) \exp\left(
\scal{x}{\psi(t,u)}\right)\;,
\end{equation}
for all $x \in D$, and for all $(t,u) \in \Rplus \times i\RR^d$ with
$ d = m+n $. We also assume that $X$ is stochastically continuous,
i.e. $X_t \to X_s$ in probability as $t \to s$.\\
Stochastic processes of this type have been studied for the first
time in the seventies, where (on the state space $D= \Rplus$) they
have been obtained as continuous-time limits of classic
Galton-Watson branching processes with and without immigration (see
\citet{Kawazu1971}). More recently, affine processes have attracted
renewed interest, due to several applications in mathematical finance, where
they are used as flexible models for asset prices, interest rates,
default intensities and other economic quantities (see
\citet{Duffie2003}
for a survey).

The process $X$ is said to be \emph{regular}, if the functions $
\Phi $ and $ \psi $ are differentiable with respect to $ t $, with
derivatives that are continuous in $ (t,u) $. This technical
condition is of crucial importance for the theory of affine
processes, as developed in \citet{Duffie2003}. Under the assumption
of regularity, it is possible to show that an affine process $X$ is
a semi-martingale (possibly killed at a state-dependent rate), and
in fact to completely characterize all affine processes in terms of
necessary and sufficient conditions on their semi-martingale
characteristics. Likewise, regularity allows to represent the
characteristic function of $X$ as the solution of so-called
generalized Riccati differential equations which determine the
fundamental functions $ \Phi $ and $ \psi $.

Without regularity it is -- a priori -- not clear how we can
``locally'' characterize the process (e.g. in terms of its
infinitesimal generator or of its semi-martingale characteristics)
and therefore such processes are not well-parameterized families
of models -- a situation which could be compared to the theory of
L\'evy processes without knowledge on the L\'evy-Khintchine
formula. Therefore in \citet{Duffie2003} the authors assume
regularity at the very beginning of their classification of affine
processes. Without the assumption of stochastic continuity there
are simple examples of non-regular Markov processes with the
affine property \eqref{Eq:affine_prop1}, based on introducing
jumps at fixed (non-random) times; see e.g.
\citet[Remark~2.11]{Duffie2003}. In contrast, if one assumes
stochastic continuity for an affine process, to the best of our
knowledge it was not now known in general whether such a process
is regular or not. Regularity has been shown in several
special cases: \citet{Kawazu1971} show automatic regularity for a
single-type continuous branching process with immigration, which
corresponds to an affine process with state space $D = \Rplus$.
\citet{Dawson2006} show regularity of affine processes under
moment conditions. In \citet{Keller-Ressel2008b} regularity of an
affine process under a mixed homogeneity and positivity condition
is shown. The approaches in the two latter publications are all
based on the general techniques from \citet{Montgomery1955} for
continuous (global) flows of homeomorphisms, which
have been extended in \citet{Filipovic2003} to continuous (local) semi-flows.\\

One of the most fascinating aspects of the regularity problem for
affine processes is the close connection to the Hilbert's fifth problem, whose
history and mathematical development can be
found in \citet{Montgomery1955}. The reason is that the functions $
\Phi $ and $ \psi $ defined by \eqref{Eq:affine_prop1} satisfy
certain functional equations, namely
\begin{equation}\label{Eq:flow1}
 \psi(t+s,u)=\psi(t,\psi(s,u)), \quad \Phi(t+s,u) =
\Phi(t,\psi(s,u)) \cdot  \Phi(s,u),\,  \text{ for } s,t \geq 0
\end{equation}
with initial conditions $ \psi(0,u) = u $ and $ \Phi(0,u)=1 $ and for all $ u$
in $\cQ$, a large enough subset of $
\CC^d $ (see Section \ref{basic_notions} for the exact definition of
$\cQ$). Such functional equations are (non-linear) relatives of the
multiplicative Cauchy functional equation
\begin{equation}\label{Eq:Cauchy1}
A(t+s)=A(t)A(s),\qquad A(0) = \operatorname{id}_d\;.
\end{equation}
formulated for instance in a set of $ d \times d $ matrices $ M_d(\CC) $, simply
by defining $ \psi(t,u)=A(t)u $. It is
well-known that the only continuous solutions of the Cauchy equation
are the exponentials $ A(t)= \exp(t \beta) $, where $ \beta $ is a $ d \times d
$ matrix.
In particular all continuous solutions of \eqref{Eq:Cauchy1} are automatically
differentiable (even analytic) with respect to $ t $, with derivatives
that are continuous in $u$.

Hilbert's fifth problem asks whether this assertion can be
extended to more general functional equations such as \eqref{Eq:flow1}:
Assuming that $\psi(t,u)$ and $\Phi(t,u)$ satisfy
\eqref{Eq:flow1}, and are differentiable in $u$, are they
necessarily differentiable in $t$? The problem has been answered
positive; for a precise formulation and a proof in the more
general context of transformation groups see
\citet[Chapter~V.5.2, e.g. Theorem~3]{Montgomery1955}. From the point of view of
stochastics this
means that moment conditions (roughly speaking the existence of a
first moment of $X$ means differentiability of $ \psi $ with
respect to $ u $) imply $t$-differentiability of $\psi$ and
$\Phi$, and thus regularity. This has already been observed in
\citet{Dawson2006}. However, moment conditions on the process $X$
are not natural in the present context of affine processes, and
the question remained open, whether the moment conditions in
\cite{Dawson2006} can simply be dropped.

We show in this article that the answer is again positive in the
general case, and that any stochastically continuous affine
process is automatically regular. There is one well-known case
where moment conditions can be dropped, namely homogeneous affine
processes, i.e., affine processes with $ \psi(t,u) = u $ for $
(t,u) \in \Rplus \times i\RR^d$. In this case $ \Phi $ simply
satisfies Cauchy's functional equation $ \Phi(t+s,u) = \Phi(t,u) \Phi(s,u) $
and therefore regularity, that is differentiability of $ \Phi $, follows by the
classical
result on \eqref{Eq:Cauchy1}. This is precisely the case of L\'evy
processes. Whence -- in the light of differentiability of
functional equations of the type \eqref{Eq:flow1} -- the results
of this paper truly extend \cite{Montgomery1955} since no
differentiability in $ u $ is assumed. We introduce the basic
definitions and some notation in Section \ref{basic_notions} and
show preliminary results in Section \ref{preliminary_results}. In
Section \ref{condition_A} we show regularity of an affine process
subject to a condition of `semi-homogeneity'. Finally we show in
Section \ref{all} that all stochastically continuous, affine
processes defined on the domain $D = \Rplus^m \times \RR^n$ are
regular by reducing the general question to regularity of
semi-homogeneous processes of Section~\ref{condition_A}.

\section{Affine processes}\label{basic_notions}
\begin{defn}[Affine process]\label{Def:affine_process} An affine process is a
time-homogeneous Markov process $(X_t,\PP^x)_{t \geq 0, x \in D}$
with state space $D = \Rplus^m \times \RR^n$, whose characteristic
function is an exponentially-affine function of the state vector.
This means that there exist functions $\Phi: \Rplus
\times i\RR^d \to \CC$ and $\psi: \Rplus \times i\RR^d \to \CC^d$
such that
\begin{equation}\label{Eq:Phi_definition}
\EE^x\left[e^{\scal{X_t}{u}}\right] = \Phi(t,u) \cdot
\exp\left(\scal{x}{\psi(t,u)}\right)\;,
\end{equation}
for all $x \in D$, and for all $(t,u) \in \Rplus \times i\RR^d$.
\end{defn}
\begin{rem}
The set $i\RR^d$ denotes the purely imaginary numbers in $\CC^d$,
that is $i\RR^d = \set{u \in \CC^d: \Re\,u = 0}$.
\end{rem}
\begin{rem}\label{rem:phi}
The above definition differs in one detail from the definition given
in \citet{Duffie2003}: In their article the right hand side of
\eqref{Eq:Phi_definition} is defined in terms of a function
$\phi(t,u)$ as $\exp\left(\phi(t,u) + \scal{x}{\psi(t,u)}\right)$,
whereas we formulate the equation in terms of $\Phi(t,u) =
\exp(\phi(t,u))$. In particular, we do not assume a priori that $ \Phi(t,u) \neq
0 $.
Their difference is subtle, but will play a role in
Lemma~\ref{Lem:phi_psi_extension} below, where we extend $\Phi(t,u)$
to a larger subset $\cQ$ of the complex numbers. Essentially, the
advantage of using $\Phi(t,u)$ is the following: if $\Phi(t,u)$ is
well-defined on a set that is not simply connected, its logarithm
$\phi(t,u)$ might only be defined as a multivalued function.
Note that our definition using $\Phi(t,u)$ is very close to that of
\citet{Kawazu1971}.
\end{rem}

\begin{ass}
We will assume throughout this article that $X$ is stochastically
continuous, i.e. for $t \to s$, the random variables $X_t$ converge
to $X_s$ in probability, with respect to all $(\PP^x)_{x \in D}$.
\end{ass}

Note that the existence of a filtered space $(\Omega, \cF,
(\cF_t)_{t \geq 0})$, where the process $\Xt$ is defined, is already
implicit in the notion of a \emph{Markov process} (we largely follow
\citet[Chapter~III]{Rogers1994} in our notation and precise
definition of a Markov process). $\PP^x$ represents the law of the
Markov process $\Xt$ \emph{started at $x$}, i.e. we have that $X_0
= x$, $\PP^x$-almost surely.

Let us at this point introduce some additional notation: We write
\[I = \set{1, \dotsc, m}\quad \text{and} \quad J = \set{m+1, \dotsc, m+n}\]
for the index sets of the $\Rplus^m$-valued component and the
$\RR^n$-valued component of $X$ respectively. For some vector $x \in
\RR^d$ we denote by $x = (x_I,x_J)$ its partition in the
corresponding subvectors, and similarly for the function $\psi(t,u)
= (\psi_I(t,u),\psi_J(t,u))$. Also, $x_i$ denotes the $i$-th element
of $x$, and $(e_i)_{i \in \set{1, \dotsc, d}}$ are the unit vectors
of $\RR^d$. We will often write
\[f_u(x) := \exp\left(\scal{u}{x}\right)\]
for the exponential function with $u \in \CC^d$ and $x \in D$. A
special role will be played by the set
\begin{equation}\label{Eq:cU_def}
\cU = \set{u \in \CC^d: \Re\,u_I \le 0, \quad \Re\,u_J = 0}\;;
\end{equation}
note that $\cU$ is precisely the set of all $u \in \CC^d$, for
which $x \mapsto f_u(x)$ is a bounded function on $D$. We also
define
\begin{equation}\label{Eq:cU_boundary_interior}
\cU^\circ = \set{u \in \CC^d: \Re\,u_I < 0, \quad \Re\,u_J = 0}.
\end{equation}


\begin{lem}\label{Lem:phi_psi_extension}
Let $\Xt$ be an affine process. Then
\begin{equation}\label{Eq:cO_def}
\cQ = \set{(t,u) \in \Rplus \times \cU: \Ex{0}{f_u(X_t)} \neq 0}\;,
\end{equation}
is open in $\Rplus \times \cU$ and there exists a unique continuous
extension of $\Phi(t,u)$ and $\psi(t,u)$ to $\cQ$, such that
\eqref{Eq:Phi_definition} holds for all $(t,u) \in \cQ$. If $(t,u)
\in (\Rplus \times \cU) \setminus \cQ$, then $\Ex{x}{f_u(X_t)} = 0$
for all $x \in D$.
\end{lem}
\begin{rem}
From the facts that $\set{0} \times \cU \in \cQ$ and that $\cQ$ is open in
$\Rplus \times \cU$ we can deduce the following: For every $u \in \cU$ there
exists $t_*(u) > 0$, such that $(t,u) \in \cQ$ for all $t \in [0,t_*(u))$.
\end{rem}
\begin{proof}
We adapt the proof of \citet[Lemma~3.1]{Duffie2003}: For $(t,u,x)
\in \Rplus \times \cU \times D$ define $g(t,u,x) =
\Ex{x}{f_u(X_t)}$. We show that for fixed $x \in D$ the function
$g(t,u,x)$ is jointly continuous in $(t,u)$: Let $(t_k,u_k)$ be a
sequence converging in $\cU$ to $(t,u)$. By stochastic continuity of
$X$ it holds that $X_{t_k} \to X_t$ in probability $\PP^x$, and thus
also in distribution. By dominated convergence we may therefore
conclude that
\[g(t_k,u_k,x) = \Ex{x}{f_{u_k}(X_{t_k})} \to \Ex{x}{f_{u}(X_{t})} =
g(t,u,x)\;,\] and thus that $g(t,u,x)$ is continuous in $(t,u)$. It
follows that $\cQ$ is open in $\Rplus \times \cU$. Because of the
affine property \eqref{Eq:Phi_definition} it holds that
\begin{equation}\label{Eq:functional_affine}
g(t,u,x) g(t,u,\xi) = g(t,u,x + \xi) g(t,u,0)
\end{equation}
for all $(t,u) \in \Rplus \times i\RR^d$ and $x,\xi \in D$. But
both sides of \eqref{Eq:functional_affine} are continuous
functions of $u \in \cU$, and moreover analytic in $\cU^\circ$.
(This follows from well-known properties of the Laplace transform
and the extension to its strip of regularity, cf.
\citet[Lemma~A.2]{Duffie2003}.) By the Schwarz reflection
principle, \eqref{Eq:functional_affine} therefore holds for all $u
\in \cU$. Assume now that $(t,u) \in (\Rplus \times \cU)
\setminus\cQ$, such that $g(t,u,0) = 0$. Then it follows from
\eqref{Eq:functional_affine} that $\Ex{x}{f_{u}(X_{t})} = g(t,u,x)
= 0$ for all $x \in D$, as claimed in the Lemma. On the other
hand, for all $(t,u) \in \cQ$ it holds that $\Phi(t,u) \neq 0$,
such that we can define $h(x) = \Phi(t,u)^{-1}g(t,u,x)$. The
function $h(x)$ is measurable and satisfies $h(x)h(\xi) = h(x +
\xi)$ for all $x,\xi \in D$. Moreover $h(0) \neq 0$ by definition
of $\cQ$. Using a standard result on measurable solutions of the
Cauchy equation (cf. \citet[Sec.~2.2]{Aczel1966}) we conclude that
there exists a unique continuous extension of $\psi(t,u)$ such
that $\Phi(t,u)^{-1}g(t,u,x) = e^{\scal{\psi(t,u)}{x}}$, and the
proof is complete.
\end{proof}
From this point on, $\Phi(t,u)$ and $\psi(t,u)$ are defined on all of $\cQ$,
and given by the unique continuous extensions of
Lemma~\ref{Lem:phi_psi_extension}. We can now give a
precise definition of regularity of an affine process:

\begin{defn}\label{Def:regularity}
An affine process $X$ is called regular, if the derivatives
\begin{equation}
F(u) =  \left.\pd{}{t} \Phi(t,u)\right|_{t = 0}, \qquad R(u) =
\left. \pd{}{t} \psi(t,u)\right|_{t = 0}
\end{equation}
exist and are continuous functions of $u \in \cU$.
\end{defn}
\begin{rem}
In \citet{Duffie2003} $F(u)$ is defined in a slightly different way,
as the derivative of $\phi(t,u)$ at $t = 0$. However the definitions
are equivalent, since $\Phi(t,u) = \exp(\phi(t,u))$ and thus
\[\left.\pd{}{t} \Phi(t,u)\right|_{t = 0} = e^{\phi(0,u)} \cdot \left.\pd{}{t}
\phi(t,u) \right|_{t = 0} = \left.\pd{}{t} \phi(t,u) \right|_{t = 0}.\]
\end{rem}

\begin{prop}\label{Prop:phi_psi_prop}
Let $ X $ be a stochastically continuous affine process. The functions $\Phi$ and
$\psi$ have the following properties:
\begin{enumerate}[(i)]
\item $\Phi$ maps $\cQ$ to the unit disc $\set{u \in \CC: |u| \le
1}$. \item $\psi$ maps $\cQ$ to $\cU$. \item $\Phi(0,u) = 1$ and
$\psi(0,u) = u$ for all $u \in \cU$. \item $\Phi$ and $\psi$ enjoy
the `\textbf{semi-flow property}': Suppose that $t,s \ge 0$ and
\linebreak[1] $(t+s,u) \in \cQ$. Then also $(t,u) \in \cQ$ and
$(s, \psi(t,u)) \in \cQ$, and it holds that
\begin{equation}\label{Eq:flow_prop}
\begin{split}
\Phi(t+s,u) &= \Phi(t,u) \cdot  \Phi(s,\psi(t,u)),\\
\psi(t+s,u) &= \psi(s,\psi(t,u)).
\end{split}
\end{equation}
\item $\Phi$ and
$\psi$ are jointly continuous on $\cQ$. \item With the remaining
arguments fixed, $u_I \mapsto \Phi(t,u)$ and $u_I \mapsto \psi(t,u)$
are analytic functions in $\set{u_I : \Re\,u_I < 0; (t,u) \in \cQ}$.
\item Let $(t,u), (t,w) \in \cQ$ with $\Re\,u \le \Re\,w$. Then
\begin{align*}
\left|\Phi(t,u)\right| &\le \Phi(t,\Re\,w)\\
\Re\,\psi(t,u) &\le \psi(t,\Re\,w)\;.
\end{align*}
\end{enumerate}
\end{prop}
\begin{proof}
Let $(t,u) \in \cQ$. Clearly $\left|\Ex{x}{f_u(X_t)}\right| \le
\Ex{x}{\left|f_u(X_t)\right|} \le 1$. On the other hand
$\Ex{x}{f_u(X_t)} = \Phi(t,u) f_{\psi(t,u)}(x)$ by
Lemma~\ref{Lem:phi_psi_extension}. Since $\norm{f_u}_\infty \le 1$
if and only if $u \in \cU$, we conclude that $\left|\Phi(t,u)\right|
\le 1$ and $\psi(t,u) \in \cU$ for all $(t,u) \in \cQ$ and have
shown (i) and (ii). Assertion (iii) follows immediately from
$\Ex{x}{f_u(X_0)} = f_u(x)$. For (iv) suppose that $(t+s,u) \in
\cQ$, such that
\begin{equation}\label{Eq:semiflow_eq1}
\Ex{x}{f_u(X_{t+s})} = \Phi(t+s,u) f_{\psi(t+s,u)}(x) \end{equation}
by Lemma~\ref{Lem:phi_psi_extension}. Applying the law of iterated
expectations and the Markov property of $X$ it holds that
\begin{equation}\label{Eq:iterated_expectation}
\Ex{x}{f_u(X_{t+s})} = \Ex{x}{\Excond{x}{f_u(X_{t+s})}{\cF_s}} =
\Ex{x}{\Ex{X_s}{f_u(X_{t})}}\;.
\end{equation}
If $(t,u) \not \in \cQ$ then the inner expectation (and consequently
the whole expression) evaluates to $0$, which is a contradiction to
the fact that $(t+s,u) \in \cQ$. It follows that
$\Ex{X_s}{f_u(X_{t})} = \Phi(t,u)f_{\psi(t,u)}(X_s)$. If
$(s,\psi(t,u)) \not \in \cQ$ then the outer expectation in
\eqref{Eq:iterated_expectation} evaluates to $0$, which is also a
contradiction. Thus also $(s,\psi(t,u)) \in \cQ$, as claimed, and we
can write \eqref{Eq:iterated_expectation} as
\begin{equation*}
\Ex{x}{f_u(X_{t+s})} = \Ex{x}{\Phi(t,u) f_{\psi(t,u)}(X_s)} =
\Phi(s,u) \cdot \Phi(t,\psi(s,u)) f_{\psi(t,\psi(s,u))}(x)\;,
\end{equation*}
for all $x \in D$. Comparing with \eqref{Eq:semiflow_eq1} the
semi-flow equations \eqref{Eq:flow_prop} follow. Assertions (v) and (vi) can be
derived
directly from the proof of Lemma~\ref{Lem:phi_psi_extension}. To
show (vii) note that
\[\left|\Ex{x}{f_u(X_t)}\right| \le \Ex{x}{\left|f_u(X_t)\right|} =
\Ex{x}{f_{(\Re\,u)}(X_t)} \le \Ex{x}{f_{(\Re\,w)}(X_t)}\;,\]
for all $x \in D$. If $(t,u)$ and $(t,w)$ are in $\cQ$, we deduce
from the affine property \eqref{Eq:Phi_definition} that
\[\left|\Phi(t,u)\right|\cdot \exp\left(\scal{x}{\Re\,\psi(t,u)}\right) \le
\Phi(t,\Re\,w) \cdot \exp\left(\scal{x}{\psi(t,\Re\,w)}\right)\;.\]
Inserting first $x = 0$ and then $C e_i$ with $C > 0$ arbitrarily
large yields the assertion.
\end{proof}

Finally we show one additional technical property concerning the
existence of derivatives of $ \Phi $ and $ \psi $ with respect to $
u $, on the interior of $ \mathcal{U} $.
\begin{lem}\label{Lem:partial_derivatives_1}Let $ X $ be a stochastically
continuous affine process. For $i \in I$, the derivatives
\[\pd{}{u_i}\Phi(t,u), \qquad \pd{}{u_i}\psi(t,u)\] exist and are continuous for
$(t,u) \in \cU^\circ \cap \cQ$.
\end{lem}

\begin{proof}
Let $i \in I$ and let $K$ be a compact subset of $\cU^\circ$. It
holds that
\begin{equation}
\left|\pd{}{u_i} \exp\left(\scal{u}{X_t}\right) \right| =
\left|X^i_t \right|\cdot \exp\left(\scal{\Re\,u}{X_t}\right)\;.
\end{equation}
The right hand side is uniformly bounded for all $t \in \Rplus, u
\in K$, and thus in particular uniformly integrable. We may conclude
that $\pd{}{u_i} \Ex{x}{e^{\scal{X_t}{u}}}$ exists and is a
continuous function of $(t,u) \in \Rplus \times K$ for any $x \in
D$. If in addition $(t,u) \in \cQ$, then
Lemma~\ref{Lem:phi_psi_extension} states that
$\Ex{x}{e^{\scal{X_t}{u}}} = \Phi(t,u)
\exp\left(\scal{x}{\psi(t,u)}\right)$. Since $K$ was an arbitrary
compact subset of $\cU^\circ$ the claim follows.
\end{proof}

\section{Affine Processes are Feller Processes}\label{preliminary_results}

In this section we prove the Feller property for all affine
processes. For \emph{regular} affine processes, this has been
shown in \citet{Duffie2003}; here we give a proof that does not
require a regularity assumption. The key to the proof are the
following properties of the function $\psi(t,u)$ for a given stochastically
continuous affine process $ X $, which will also
be used in the proof of our main result in Section~\ref{all}:
\begin{description}\label{properties_a_b}
\item[Property A] $\psi(t,.)$ maps $\cU^\circ$ to $\cU^\circ$.
\item[Property B] $\psi_J(t,u) = e^{\beta t}u_J$ for all $(t,u) \in \cQ$, with
$\beta$ a real $n \times n$-matrix.
\end{description}
Let us give here an intuitive example illustrating the second
property, which has already been observed in \citet[Prop~2.1 and
Cor.~2.1]{Dawson2006}: Consider an affine process with one-dimensional state
space $D = \RR$, and the property that $\Phi(t,u) = 1$. Then for any
initial value $x \in \RR$
\[\Ex{x}{f_{iy}(X_t)} = e^{x\psi(t,iy)} \quad \text{and} \quad
\Ex{-x}{f_{iy}(X_t)} = e^{-x \psi(t,iy)}\]
are both characteristic functions, and moreover reciprocal to each
other. But a well-known result (cf.~\citet[Thm.~2.1.4]{Lukacs1960}) states that
the only
characteristic functions, whose reciprocals are also
characteristic functions correspond to degenerate distributions, i.e.~Dirac
measures. Here, this implies that $\psi(t,iy) = iy m(t)$, for $m(t)$ a
deterministic function. Moreover, by the Markov property $m(t+s) =
m(t) m(s)$, which is Cauchy's functional equation with the unique
continuous solution $m(t) = e^{\lambda t}m(0)$, for some $ \lambda \in
\mathbb{R} $. Hence,
$\psi(t,iy)$ is necessarily of the form $e^{\lambda t}\,iy$ and
satisfies therefore \emph{Property~B}.

As we shall show, the argument can be extended to the case of arbitrary
$\Phi(t,u)$ and
to the general state space $D = \Rplus^m \times \RR^n$. The next
Lemma is the first step in this direction:

\begin{lem}\label{Lem:linearity}
Let $ X $ be a stochastically continuous affine process. Let $K \subseteq
\set{1, \dotsc, d}$, $k \in \set{1, \dotsc, d}$,
and let $(t_n)_{n \in \NN}$ be a sequence such that $t_n \downarrow
0$. Define $\Omega_K := \set{y \in \RR^d: y_i = 0\; \text{for}\;i
\not \in K}$, and suppose that
\[\Re\,\psi_k(t_n,iy) = 0 \qquad \text{for all $y \in \Omega_K$ and $n \in
\NN$}\;.\] Then there exist $\zeta_k(t_n) \in \RR^{|K|}$ and an
increasing sequence of positive numbers $R_n$ such that $R_n
\uparrow \infty$ and
\[\psi_k(t_n,iy) = \scal{\zeta_k(t_n)}{i y_K}\;\]
for all $y \in \Omega_K$ with $|y| < R_n$.
\end{lem}

For the proof we will use the following result:
\begin{lem}\label{Lem:charf}
Let $\Theta$ be a positive definite function on $\RR^d$ with
$\Theta(0) = 1$. Then
\[\left|\Theta(y+z) - \Theta(y)\Theta(z)\right|^2 \le \left(1 -
|\Theta(y)|^2\right)\left(1 - |\Theta(z)|^2\right) \le 1\]
for all $y,z \in \RR^d$.
\end{lem}
\begin{proof}[Proof of Lemma~\ref{Lem:charf}]
The result follows from considering the matrix
\[M_\Theta(y,z) := \left(
\begin{array}{ccc}
  \Theta(0) & \overline{\Theta(y)} & \Theta(z) \\
  \Theta(y) & \Theta(0) & \Theta(y+z) \\
  \overline{\Theta(z)} & \overline{\Theta(y+z)} & \Theta(0) \\
\end{array}
\right), \qquad y,z \in \RR^d, y \neq z\;,\] which is positive
semi-definite by definition of $\Theta$. The inequality is then
derived from the fact that $\det M_\Theta(y,z) \ge 0$. See
\citet[Lemma~3.5.10]{Jacob2001} for details.
\end{proof}

\begin{proof}[Proof of Lemma~\ref{Lem:linearity}]
As the characteristic function of the (possibly defective) random
variable $X_{t_n}$ under $\PP^x$, the function $y \mapsto \Ex{x}{
f_{iy}(X_{t_n})}$ is positive definite for any $x \in D, n \in \NN$.
We define now for every $y \in \Omega_K$, $c > 0$, and $n \in \NN$,
the function
\[\Theta(y;n,c) := \frac{1}{\Phi(t_n,0)} \Ex{c e_k}{f_{iy}(X_{t_n})} =
\frac{\Phi(t_n,iy)}{\Phi(t_n,0)} \exp\Big(c \cdot
\psi_k(t_n,iy)\Big)\;.\] 
Clearly, as a function of $y \in \Omega_K$,
also $\Theta(y;n,c)$ is positive definite. In addition it satisfies
$\Theta(0;n,c) = \exp\left(c \cdot \psi_k(t_n,0)\right) = 1$, for large enough
$n$, say $n \ge N$, by the following argument: It should be obvious, that
$\exp(c \cdot \psi_k(t_n,0)$ is always a real quantity. By assumption,
$\psi_k(t_n,0)$ is purely imaginary, such that it must be an integer multiple of
$\pi$ for all $n \in \NN$. But $\psi_k(0,0) = 0$, and $\psi_k(t,0)$ is
continuous in $t$ by Proposition~\ref{Prop:phi_psi_prop}, and we conclude
$\psi_k(t_n,0) = 0$ for large enough $n$. 

Thus, for $n \ge N$, we may apply Lemma~\ref{Lem:charf} to $\Theta$, and
conclude that for any
$y,z \in \Omega_K$, $c > 0$ and $n \ge N$
\begin{equation}\label{Eq:psi_interm_1}
\left|\Theta(y + z;n,c) - \Theta(y;n,c) \cdot \Theta(z;n,c)\right|^2
\le 1\;.
\end{equation}
For compact notation we define the abbreviations
\begin{align*}
&r_1 = \left|\frac{\Phi(t_n,i(y+z))}{\Phi(t_n,0)}\right|, &\qquad
&r_2 =
\left|\frac{\Phi(t_n,iy)\Phi(t_n,iz)}{\Phi(t_n,0)^2}\right|\,,\\
&\alpha_1 = \mathrm{Arg}\,\frac{\Phi(t_n,i(y+z))}{\Phi(t_n,0)},
&\qquad &\alpha_2 =
\mathrm{Arg}\,\frac{\Phi(t_n,iy)\Phi(t_n,iz)}{\Phi(t_n,0)^2}\,,\\
&\beta_1 = \Im\,\psi_k(t_n,i(y+z)), &\qquad &\beta_2 = \Im\,\psi_k(t_n,iy) +
\Im\,\psi_k(t_n,iz),
\end{align*}
where we suppress the dependency on $y,t,z$ for the moment. It
holds that
\[\left|r_1 e^{(\alpha_1 + c \beta_1)i} -
r_2 e^{(\alpha_2 + c \beta_2)i}\right|^2 = r_1^2 + r_2^2 - 2r_1 r_2
\cos\Big(\alpha_1 - \alpha_2 + (\beta_1 -
\beta_2)c\Big)\;.\]
Using the elementary inequality $2 r_1 r_2 \le r_1^2 + r_2^2$, we derive
\begin{equation*}\label{Eq:complex_eq}
2 r_1 r_2 \left\{1 - \cos\Big(\alpha_1 - \alpha_2 + (\beta_1 -
\beta_2)c\Big) \right\} \le \left|r_1 e^{(\alpha_1 + c \beta_1)i} -
r_2 e^{(\alpha_2 + c \beta_2)i}\right|^2\;,
\end{equation*}
which combined with inequality \eqref{Eq:psi_interm_1} yields
\begin{equation}\label{Eq:link}
r_1 r_2 \left(1 - \cos\Big(\alpha_1 - \alpha_2 + (\beta_1 -
\beta_2)c\Big)\right) \le \frac{1}{2}\;.
\end{equation}
Define now $R_n = 0$ for $n < N$, and 
\[R_n := \sup \set{\rho \ge 0: r_1(y,t_n,z)r_2(y,t_n,z) > \frac{1}{2} \;
\text{for}\;y,z \in \Omega_K\;\text{with}\;|y| \le \rho, |z|  \le \rho}\]
for $n \ge N$. Note that $R_n \uparrow \infty$: This follows from
$r_1(y,0,z) = r_2(y,0,z) = 1$ for all
$y,z \in \Omega_K$, and the continuity of $r_1$ and $r_2$.

Suppose that
\[\beta_1 - \beta_2 = \Im\,\psi_k(t_n,i(y+z)) - \Im\,\psi_k(t_n,iy) -
\Im\,\psi_k(t_n,iz) \neq 0\]
for any $n \in \NN$ and $y,z \in \Omega_K$ with $|y| < R_n$, $|z| <
R_n$. Then there exists an $c
> 0$ such that
\[\cos \Big(\alpha_1 - \alpha_2 + (\beta_1 -
\beta_2)c\Big) = -1\,.\]Inserting into \eqref{Eq:link} we obtain
\[\frac{1}{2}\cdot 2 < r_1 r_2 \left(1 - \cos \Big(\alpha_1 - \alpha_2 +
(\beta_1 -
\beta_2)c\Big)\right) \le \frac{1}{2}\;,\] a contradiction. We
conclude that
\begin{equation}\label{Eq:Cauchy_eq}
\beta_1 - \beta_2 = \Im\,\psi_k(t_n,i(y+z)) - \Im\,\psi_k(t_n,iy) -
\Im\,\psi_k(t_n,iz) = 0\;,
\end{equation}
for all $y,z \in \Omega_K$ with $|y| < R_n$, $|z| < R_n$. Equation
\eqref{Eq:Cauchy_eq} is nothing but Cauchy's first functional
equation. Since $\psi(t,.)$ is continuous, it follows that
$\Im\,\psi_k$ is a linear function of $y_K$. In addition, $\Re\,\psi_k(t_n,y)$
is zero, by assumption, such that there exists some
real vector $\zeta_k(t_n)$ with 
\begin{equation}\label{Eq:linear_interm}
\psi_k(t_n,iy) = \scal{\zeta_k(t_n)}{ i y_K}\;.\end{equation} for
all $y \in \Omega_K$ with $|y| < R_n$, and the
Lemma is proved.
\end{proof}

We use the above Lemma to show the following Proposition, which
implies Property~B of $\psi$, that was introduced at the beginning
of the section:

\begin{prop}\label{Prop:J_exp}
Let $\Xt$ be a stochastically continuous affine process on $D = \Rplus^m \times
\RR^n$ and
denote by $J$ its real-valued components. Then there exists a real
$n \times n$-matrix $\beta$ such that $\psi_J(t,u) = e^{t \beta}
u_J$ for all $(t,u) \in \cQ$.
\end{prop}
\begin{proof}
Consider the definition of $\cU$ in \eqref{Eq:cU_def}. Since
$\psi(t,u)$ takes by Proposition~\ref{Prop:phi_psi_prop} values in
$\cU$ it is clear that $\Re\,\psi_J(t,iy) = 0$ for \emph{any} $(t,y)
\in \Rplus \times \RR^d$. Fix now some $t_* > 0$ and define $t_n :=
t_*/n$ for all $n \in \NN$. We can apply Lemma~\ref{Lem:linearity}
with $K = \set{1, \dotsc, d}$ and any choice of $k \in J$, to obtain
a sequence $R_n \uparrow \infty$ (even independent of $ k $), such
that
\begin{equation}\label{Eq:psiJ_linear}
\psi_J(t_n,iy) = \Xi(t_n) \cdot iy\;,
\end{equation}
for all $y \in \RR^d$ with $|y| < R_n$. Here $\Xi(t_n)$ denotes the
real $n \times d$-matrix formed by concatenating the column vectors
$(\zeta_k(t_n))_{k =1, \dotsc, d}$ obtained
from Lemma~\ref{Lem:linearity}.

Let $i \in I$, $n \in \NN$, define $\Omega_n := \set{\omega \in \CC:
|\omega| \le R_n , (t_n,e_i \omega) \in \cQ}$, and consider the
function
\[h_n:  \Omega_n \to \CC^n: \quad \omega \mapsto \psi_J(t_n,\omega
e_i) - \Xi(t) \cdot \omega e_i\;.\]

This is an analytic function on $\Omega_n^\circ$ and continuous on
$\Omega_n$. According to the Schwarz reflection principle, $h_n$ can
be extended to an analytic function on an open superset of
$\Omega_n$. But \eqref{Eq:psiJ_linear} implies that the function
$h_n$ takes the value $0$ on a subset with an accumulation point in
$\CC$. We conclude that $h_n$ is zero everywhere. In particular we
have that
\[0 = \Re\,\psi_J(t_n,\omega e_i) -  \Xi(t_n) \cdot \Re\,\omega e_i = \Xi(t_n)
\cdot \Re\,\omega e_i \;,\]
for all $\omega \in \Omega_n$. This can only hold true, if the
$i$-th column of $\Xi(t_n)$ is zero. Since $i \in I$ was arbitrary
we have reduced \eqref{Eq:psiJ_linear}
to\begin{equation}\label{Eq:psiJ_linear_2} \psi_J(t_n,u) =
\Xi_0(t_n) \cdot u_J\;,
\end{equation}
for all $(t_n,u) \in \cQ$, such that $|u_J| < R_n$. Here
$\Xi_0(t_n)$ denotes the $n \times n$-submatrix of
$\Xi(t_n)$ that results from dropping the zero-columns.

Fix an arbitrary $u_* \in \cU$ with $(t_*,u_*) \in \cQ$. By
Proposition \ref{Prop:phi_psi_prop} we know that also $(t,u_*) \in
\cQ$ for all $t \in [0,t_*]$, such that \linebreak[1] $R := \sup
\set{|\psi_K(t,u_*)|: t \in [0,t_*]}$ is well-defined. Since
$\psi(t,u)$ is continuous, $R$ is finite. Choose $N$ such that $R_n
> R$ for all $n \ge N$. Using the semi-flow equation we can write
$\psi_J(t_*,u_*)$ as
\begin{multline}\label{Eq:self_insertion}
\psi_J(t_*,u_*) = \psi_J\left(t_n
,\psi(t_* \tfrac{n-1}{n},u_*)\right) = \\
= \Xi_0(t_n) \cdot \psi_J(t_* \tfrac{n-1}{n}, u_*) = \dotsm =
\Xi_0(t_n)^n \cdot u_*\;;
\end{multline}
for any $n \ge N$. Thus, the functional equation $\psi(t,u) =
\Xi_0(t) \cdot u_J$ actually holds for all $(t,u) \in \cQ$. Another
application of the semi-flow property yields then, that
\[\Xi_0(t+s) = \Xi_0(t) \Xi_0(s), \qquad \text{for all $t,s \geq 0$}\;.\]
Since $\Xi_0(0) = 1$, $\Xi_0$ is continuous and satisfies the second
Cauchy functional equation, it follows that $\Xi_0(t) = e^{\beta t}$
for some real $n \times n$-matrix $\beta$, which completes the proof.
\end{proof}

The next proposition shows that also Property~A holds true for
$\psi$, as we have claimed at the beginning of the section.

\begin{prop}\label{Prop:interior_to_interior}
Suppose that $(t,u) \in \cQ$. If $u \in \cU^\circ$, then $\psi(t,u)
\in \cU^\circ$.
\end{prop}
\begin{proof}
For a contradiction, assume there exists $(t,u) \in \cQ$ such that
$u \in \cU^\circ$, but $\psi(t,u) \not \in \cU^\circ$. This implies
that there exists $k \in I$, such that $\Re\,\psi_k(t,u) = 0$. Let
$\cQ_{t,k} = \set{\omega \in \CC: \Re\,\omega \le 0;\, (t,\omega
e_k) \in \cQ}$. From the inequalities of
Proposition~\ref{Prop:phi_psi_prop}.vii we deduce that
\begin{equation}\label{Eq:sandwich}
0 = \Re\,\psi_k(t,u) \le \psi_k(t,\Re\,\omega \cdot e_k) \le 0\;,
\end{equation}
and thus that $\psi_k(t,\Re\,\omega \cdot e_k) = 0$ for all $\omega
\in \cQ_{t,k}$ with $\Re\,u_k \le \Re\,\omega$. By
Proposition~\ref{Prop:phi_psi_prop}.(vi), $\psi_k(t,\omega e_k)$ is
an analytic function of $\omega$. Since it takes the value zero on a
set with an accumulation point, it is zero everywhere, i.e.
$\psi_k(t,\omega e_k) = 0$ for all $\omega \in \cQ_{t,k}$. The same statement
holds true for $t$ replaced by $t/2$: Set
$\lambda := \Re\,\psi_k(t/2,u)$. If $\lambda = 0$, we can proceed
exactly as above, only with $t/2$ instead of $t$. If $\lambda < 0$,
then we have, by another application of
Proposition~\ref{Prop:phi_psi_prop}.vii, that
\begin{equation*}
0 = \Re\,\psi_k(t,u) = \Re\,\psi_k(t/2,\psi(t/2,u)) \le
\psi_k(t/2,\lambda e_k) \le \psi_k(t/2,\Re\,\omega e_k) \le 0\;,
\end{equation*}
for all $\omega \in \cQ_{t/2,k}$ such that $\lambda \le
\Re\,\omega$. Again we use that an analytic function that takes
the value zero on a set with accumulation point, is zero everywhere,
and obtain that $\psi_k(t/2,\omega e_k) = 0$ for all $\omega \in
\cQ_{t/2,k}$. Repeating this argument, we finally obtain a sequence
$t_n \downarrow 0$, such that
\begin{equation}\label{Eq:real_part_zero}
\psi_k(t_n,\omega e_k) = 0 \qquad \text{for all} \; \omega \in
\cQ_{t_n,k}.
\end{equation}
We can now apply Lemma~\ref{Lem:linearity} with $K = \set{k}$, which
implies that $\psi_k$ is of the linear form
\[\psi_k(t_n,\omega e_k) = \zeta_k(t_n) \cdot \omega, \quad \text{for all} \quad
\omega \in \cQ_{t_n,k}\,\text{with}\; |\omega| \le R_n,\] where
$\zeta_k(t_n)$ are real numbers, and $R_n \uparrow \infty$. Note
that since $\zeta_k(t_n) \to 1$ as $t_n \to 0$, we have that
$\zeta_k(t_n) > 0$ for $n$ large enough. Choosing now some
$\omega_*$ with $\Re\,\omega_* < 0$ it follows that
$\Re\,\psi_k(t_n,\omega_* e_k) < 0$ --  with strict inequality. This
is a contradiction to \eqref{Eq:real_part_zero}, and the assertion
is shown.
\end{proof}

We are now prepared to show the main result of this section:
\begin{thm}\label{Thm:Feller}
Every stochastically continuous affine process $ X $ is a Feller
process.
\end{thm}
\begin{rem}
As an immediate Corollary to this theorem, every stochastically continuous
affine process has a c\`adl\`ag version, see for instance \cite{Rogers1994}.
\end{rem}
\begin{proof}
By stochastic continuity of $\Xt$ and dominated convergence, it
follows immediately that $P_t f(x) = \Ex{x}{f(X_t)} \to f(x)$ as $t
\to 0$ for all $f \in C_0(D)$ and $x \in D$. To prove the Feller
property of $\Xt$ it remains to show that $P_t (C_0(D)) \subseteq
C_0(D)$: For $u_I \in \CC^m$ with $\Re\,u_I < 0$ and $g \in
C^\infty_c(\RR^n)$, i.e.~a smooth function with bounded support, define the
functions
\[h(x; u_I,g) =
e^{\scal{u_I}{x_I}}\int_{\RR^n}{f_{iy}(x_J)g(y)\,dy}\] mapping $D$
to $\CC$, and the set
\begin{equation*}
P := \set{h(x; u_I,g) \,:\; u_I \in \CC^d, \, \Re\, u_I < 0, g \in
C^\infty_c(\RR^n)}\;.
\end{equation*}
Denote by $\cL(P)$ the set of (complex) linear combinations of
functions in $P$. From the Riemann-Lebesgue Lemma it follows that
$\int_{\RR^n}{f_{iy}(x_J)g(y)\,dy}$ vanishes at infinity, and thus
that $\cL(P) \subset C_0(D)$. It is easy to see that $\cL(P)$ is a
subalgebra of $C_0(D)$, that is in addition closed under complex
conjugation and multiplication. (Note that the product of two Fourier transforms
of compactly supported functions $ g_1, g_2 $ is the Fourier transform of a
compactly supported function, namely $ g_1 * g_2 $.) It is also straight-forward
to check that $\cL(P)$ is
point separating and vanishes nowhere (i.e. there is no $x_0 \in D$
such that $h(x_0) = 0$ for all $h \in \cL(P)$).
Using a suitable version of the Stone-Weierstrass theorem
(e.g.~\citet[Corollary~7.3.9]{Semadeni1971}), it follows that
$\cL(P)$ is dense in $C_0(D)$.

Fix some $t \in \Rplus$ and let $h(x) \in P$. By
Lemma~\ref{Lem:phi_psi_extension} it holds that $\Ex{x}{f_u(X_t)} =
\Phi(t,u) \exp\left(\scal{x}{\psi(t,u)}\right)$ whenever $(t,u) \in
\cQ$, and $\EE^x[f_{u}(X_t)] = 0$ whenever $(t,u) \not \in
\cQ$. Moreover, by Proposition~\ref{Prop:J_exp} we know that
$\psi_J(t,u) = e^{\beta t}u_J$ for all $(t,u) \in \cQ$. Thus, writing $u =
(u_I,iy)$, we have
\begin{align}\label{Eq:Pth_fourier}
P_t h(x) &=
\EE^x\left[\int_{\RR^n}{f_{(u_I,iy)}(X_t)g(y)\,dy}\right] =
\int_{\RR^n}{\Ex{x}{f_{(u_I,iy)}(X_t)}g(y)\,dy} = \\
&= \int_{\set{u \in \cU: (t,u) \in \cQ}}{\Ex{x}{f_{(u_I,iy)}(X_t)} g(y)\,dy} =
\notag\\
&= \int_{\set{u \in \cU: (t,u) \in \cQ}}\Phi(t,u_I,iy) \exp\left(
\scal{x_I}{\psi_I(t,u_I,iy)} + \scal{x_J}{e^{t\beta}
iy}\right)g(y)\,dy\;.  \notag
\end{align}
Since $(u_I,iy) \in \cU^\circ$ it follows by
Proposition~\ref{Prop:interior_to_interior} that also
$\Re\,\psi_I(t,u_I,iy) < 0$ for any $y \in \RR^n$. This shows that
$P_t h(x) \to 0$ as $|x_I| \to \infty$. In addition, as a function
of $x_J$, \eqref{Eq:Pth_fourier} can be interpreted as the Fourier
transformation of a compactly supported density. The
Riemann-Lebesgue Lemma then implies that $P_t h(x) \to 0$ as $|x_J|
\to \infty$, and we conclude that $P_t h \in C_0(D)$. The assertion
extends by linearity to every $h \in \cL(P)$, and finally by the
density of $\cL(P)$ to every $h \in C_0(D)$. This proves that the
semi-group $\Pt$ maps $C_0(D)$ into $C_0(D)$, and hence that $\Xt$
is a Feller process.
\end{proof}

\section{All semi-homogeneous affine processes are regular}\label{condition_A}

\begin{defn}\label{Def:condA}
We say that a stochastically continuous affine process is
\emph{semi-homogeneous}, if for all
$x \in D$, $(t,u) \in \Rplus \times i\RR^d$
\begin{equation}\label{Eq:semihomogeneous}
\Ex{x}{e^{\scal{X_t}{u}}} = e^{\scal{x_J}{u_J}} \cdot
\Ex{(x_I,0)}{e^{\scal{X_t}{u}}}\,.
\end{equation}
\end{defn}

The above condition is equivalent to the statement that for any $y$
of the form $y = (0,y_J)$, the law of $X_t + y$ under $\PP^x$ equals
the law of $X_t$ under $\PP^{(x + y)}$. If this held true for any $y
\in D$ we would speak of a (space-)homogeneous Markov process. Since
we impose the condition only for $y$ of the form $(0,y_J)$, we call
the process semi-homogeneous. Note that semi-homogeneous affine
processes are frequently encountered in mathematical finance: Affine
stochastic volatility models (e.g. the Heston model) are typically
based on a semi-homogeneous affine process; see \citet{Keller-Ressel2008a}.

Combining Definition~\ref{Def:condA} with the affine property
\eqref{Eq:Phi_definition}, it is easy to see that the following
holds:
\begin{lem}\label{Lem:semi_homo}
A stochastically continuous affine process is semi-homogeneous, if and only if
$\psi_J(t,u) =
u_J$ for all $(t,u) \in \Rplus \times i\RR^d$ (or equivalently for
all $(t,u) \in \cQ$).
\end{lem}

The main result of this section is the following:

\begin{thm}\label{Thm:condA_regular}
Every semi-homogeneous, stochastically continuous affine process is
regular.
\end{thm}
Our proof uses the techniques originally presented in \citet{Montgomery1955} for
continuous transformation groups, and follows in part the proof of
\cite[Theorem~4.1]{Dawson2006}.
\begin{proof}To simplify calculations we embed $\Phi(t,u)$ and $\psi(t,u)$ into
the extended
semi-flow $\Upsilon(t,u)$, that is we set $\wh{\cQ} := \cQ \times
\CC$ and define
\begin{equation}\label{Eq:bigflow_1}
\Upsilon: \; \wh{\cQ} \to \CC^{d+1}, \quad (t,u_1, \dotsc
u_d,u_{d+1}) \mapsto \left(\begin{array}{c}
  \psi(t,(u_1,\dotsc,u_d)) \\
  \Phi(t,(u_1,\dotsc,u_d)) \cdot u_{d+1}
\end{array}%
\right).
\end{equation}
Note that all vectors $u$ have now a $(d+1)$-th component added;
this component will be assigned to the non-negative components $I$,
such that under slight abuse of notation we now write $I  = \set{1,
\dotsc, m, d+1}$. The semi-flow property is preserved by
$\Upsilon(t,u)$, i.e. $\Upsilon(t + s, u) =
\Upsilon(t,\Upsilon(s,u))$ for all $(t+s,u) \in \wh{\cQ}$. The
semi-homogeneity condition on $X$ implies that $\Upsilon_J(t,u) =
u_J$ for all $(t,u) \in \wh{\cQ}$. Clearly, the time derivative
exists and vanishes, i.e.,
$$
\left.\pd{}{t}\Upsilon_J(t,u)\right|_{t = 0} =0
$$
for all $u \in \cU \times \CC$. In the rest of the proof we thus
focus on the remaining (not space-homogeneous) components. Let $u
\in \wh{\cU}^\circ := \cU^\circ \times \CC$ be fixed and assume that
$t, s \in \Rplus$ are small enough such that $\Phi(t+s,u)$,
$\psi(t+s,u)$ and their $u$-derivatives are always well-defined (cf.
Lemma~\ref{Lem:partial_derivatives_1}). Denote by
$\pd{\Upsilon_I}{u_I}(t,u)$ the Jacobian of $\Upsilon_I$ with
respect to $u_I$. Using a Taylor expansion we have that
\begin{multline}\label{Eq:Taylor_flow}
\int_0^s \Upsilon_I(r,\Upsilon(t,u))\,dr - \int_0^s
\Upsilon_I(t,u)\,dr = \int_0^s \pd{\Upsilon_I}{u_I}(r,u)\,dr \cdot
\left(\Upsilon_I(t,u) - u_I\right) + \\
 + o\left(\norm{\Upsilon_I(t,u) - u_I}\right)\;.
\end{multline}
On the other hand, using the semi-flow property of $\Upsilon$ we
can write the left side of \eqref{Eq:Taylor_flow} as
\begin{equation}\label{Eq:Flow_interm1}
\begin{split}
 \int_0^s{\Upsilon_I(r,\Upsilon(t,u))\,dr} -
\int_0^s{\Upsilon_I(r,u)\,dr} &=
\int_0^s{\Upsilon_I(r + t,u)\,dr} - \int_0^s{\Upsilon_I(r,u)\,dr} = \\
= \int_t^{s + t}{\Upsilon_I(r,u)\,dr} -
\int_0^s{\Upsilon_I(r,u)\,dr} &=
 \int_s^{s + t}{\Upsilon_I(r,u)\,dr} - \int_0^t{\Upsilon_I(r,u)\,dr} =  \\
&= \int_0^t{\Upsilon_I(r + s,u)\,dr} -
\int_0^t{\Upsilon_I(r,u)\,dr}\;.
\end{split}
\end{equation}
Denoting the last expression by $I(s,t)$ and combining
\eqref{Eq:Taylor_flow} with \eqref{Eq:Flow_interm1} we obtain
\[\lim_{t \downarrow 0} \frac{\norm{\frac{1}{s}I(s,t)}}{\norm{\Upsilon_I(t,u) -
u_I}} = \norm{\frac{1}{s}\int_0^s
\pd{\Upsilon_I}{u_I}(r,u)\,dr}\;.\] Define $M(s,u) :=
\frac{1}{s}\int_0^s \pd{\Upsilon_I}{u_I}(r,u)\,dr$. Note that as $s
\to 0$, it holds that $M(s,u) \to \pd{\Upsilon_I}{u_I}(0,u) = I_I$
(the identity matrix). Thus for $s$ small enough $\norm{M(s,u)} \neq
0$, and we conclude that
\begin{multline} \label{Eq:Upsilon_lim}
\lim_{t \downarrow 0} \frac{1}{t} \norm{\Upsilon_I(t,u) - u_I} = \\
= \norm{\lim_{t \downarrow} \frac{I(s,t)}{st}} \cdot
\norm{M(s,u)}^{-1} =  \norm{\frac{\Upsilon_I(s,u) - u_I}{s}} \cdot
\norm{M(s,u)}^{-1}\;.
\end{multline}
The right hand side of \eqref{Eq:Upsilon_lim} is well-defined and
finite, implying that also the limit on the left hand side is.
Thus, combining \eqref{Eq:Taylor_flow} and
\eqref{Eq:Flow_interm1}, dividing by $st$ and taking the limit $t
\downarrow 0$ we obtain
\[\lim_{t \downarrow 0} \frac{\Upsilon_I(t,u) - u_I}{t} = \frac{\Upsilon_I(s,u)
- u_I}{s} \cdot M(s,u)^{-1}\;.\]
Again we may choose $s$ small enough, such that $M(s,u)$ is
invertible, and the right hand side of the above expression is
well-defined. The existence and finiteness of the right hand side
then implies the existence of the limit on the left. In addition the
right hand side is a continuous function of $u \in \wh{\cU}^\circ$,
such that also the left hand side is. Adding back the components
$J$, for which a time derivative trivially exists (recall that $\psi_J(t,u) =
u_J$ for all $t \ge 0$), we obtain that
\begin{equation}\label{Eq:R_limit} \cR(u) := \lim_{t \downarrow 0}
\frac{\Upsilon(t,u) - u}{t} =
\left.\pd{}{t}\Upsilon(t,u)\right|_{t = 0}
\end{equation}
exists and is a continuous function of $u \in \wh{\cU}^\circ$.
Denoting the first $d$ components of $\cR(u)$ by $R(u)$ and the
$d+1$-th component by $F(u)$ we can `disentangle' the extended
semi-flow $\Upsilon$, drop the $(d+1)$-th component of $u$, and
see that
\begin{equation}\label{Eq:FR_limit}
\left.F(u) := \pd{}{t}\Phi(t,u)\right|_{t = 0} \quad \text{and}
\quad \left.R(u) := \pd{}{t}\psi(t,u)\right|_{t = 0}
\end{equation}
are likewise well-defined and continuous on $\cU^\circ$.

To show that $\Xt$ is regular affine it
remains to show that \eqref{Eq:FR_limit} extends continuously to
$\cU$: To this end let $t_n \downarrow 0$, $x \in D$, $u \in \cU^\circ $,
and rewrite \eqref{Eq:FR_limit} as
\begin{multline}\label{Eq:FR_representasion}
F(u) + \scal{x}{R(u)} = \lim_{n \to \infty}
\frac{\Phi(t_n,u)\exp\left( \scal{x}{\psi(t_n,u) - u}\right) -
1}{t_n} =
\\ = \lim_{n \to \infty} \frac{f_{-u}(x)\,\Ex{x}{f_u(X_{t_n})}
- 1}{t_n}
 = \lim_{n \to \infty} \frac{1}{t_n} \left\{\int_{D}
e^{\scal{\xi - x}{u}}\,p_{t_n}(x,d\xi) - 1\right\} = \\
= \lim_{n \to \infty} \frac{1}{t_n} \left\{\int_{D - x}
\left(e^{\scal{\xi}{u}} - 1\right)\,\wt{p}_{t_n}(x,d\xi) +
\frac{p_{t_n}(x,D) - 1}{t_n}\right\}\;,
\end{multline}
where $p_t(x,d\xi)$ is the transition kernel of the Markov process
$\Xt$, and $\wt{p}_t(x,d\xi)$ is its `shifted transition kernel'
$\wt{p}_t(x,d\xi) := p_t(x,d\xi + x)$. The right hand side of
\eqref{Eq:FR_representasion} can be regarded as a limit of
log-characteristic functions of (infinitely divisible)
sub-stochastic measures\footnote{Note that
$\exp\left(\frac{1}{t_n}\left\{\int_{D - x} \left(e^{\scal{\xi}{u}}
- 1\right)\,\wt{p}_t(x,d\xi)\right\}\right)$ is the characteristic
function of a compound Poisson distribution with intensity
$\frac{1}{t_n}$ and jump measure $\wt{p}_t(x,d\xi)$.}. That is, there
exist infinitely divisible sub-stochastic measures $\mu_n(x,d\xi)$,
such that
\[\exp\left(F(u) + \scal{x}{R(u)}\right) = \lim_{n \to \infty}
\int_{\RR^d}{e^{\scal{u}{\xi}}\,\mu_n(x,d\xi)}, \quad \text{for all $u \in
\cU^\circ$.}\]
Let now $\theta \in \RR^d$ with $\theta_I < 0$ and $\theta_J = 0$
(note that $\theta \in \cU^\circ$) and consider the exponentially
tilted measures $e^{\scal{\theta}{\xi}}\,\mu_n(x,d\xi)$. Their
characteristic functions converge to $\exp\left(F(u + \theta) +
\scal{x}{R(u+\theta)}\right)$. Thus, by L\'{e}vy's continuity
theorem, there exists $\mu_*(x,d\xi)$ such that
$e^{\scal{\theta}{\xi}}\,\mu_n(x,d\xi) \to \mu_*(x,d\xi)$ weakly. On
the other hand, by Helly's selection theorem, $\mu_n(x,d\xi)$ has a
vaguely convergent subsequence, which converges to some measure
$\mu(x,d\xi)$. By uniqueness of the weak limit we conclude that
$\mu(x,d\xi) = e^{\scal{-\theta}{\xi}}\mu_*(x,d\xi)$. Thus we have
that for all $x \in D$ and $u \in \cU^\circ$ with $\Re\,u$ in a
neighborhood of $\theta$,
\begin{multline}\label{Eq:tilted_conv}
\exp\left(F(u) + \scal{x}{R(u)}\right) = \lim_{n \to \infty}
\int_{\RR^d}{e^{\scal{u}{\xi}}\,\mu_n(d,d\xi)} = \\
 = \lim_{n \to
\infty} \int_{\RR^d}{e^{\scal{u -
\theta}{\xi}}\,e^{\scal{\theta}{\xi}}\mu_n(x,d\xi)} =
\int_{\RR^d}{e^{\scal{u - \theta}{\xi}}\,\mu_*(x,d\xi)} =
\int_{\RR^d}e^{\scal{u}{\xi}}\,\mu(x,d\xi)\,.
\end{multline}
But the choice of $\theta$ was arbitrary, such that
\eqref{Eq:tilted_conv} extends to all $u \in \cU^\circ$. Applying
dominated convergence to the last term of \eqref{Eq:tilted_conv}
shows that both $F$ and $R$ have a continuous extension to all of
$\cU$, which we also denote by $F$ and $R$ respectively. 

It remains to show that \eqref{Eq:FR_limit}
remains valid on $\cU$: Let $u \in \cU$ and $(u_n)_{n \in \NN} \in \cU^\circ$
such that
$u_n \to u$. Remember that by
Proposition~\ref{Prop:interior_to_interior} $u_n \in \cU^\circ$
implies that also $\psi(t,u_n) \in \cU^\circ$ for any $t \ge 0$.
Thus we have
\begin{multline}\label{Eq:R_ext_is_der}
\int_0^t{R(\psi(s,u))\,ds} = \int_0^t{\lim_{u_n \to
u}R(\psi(s,u_n))\,ds} = \lim_{u_n \to u}
\int_0^t{R(\psi(s,u_n))\,ds} = \\
= \lim_{u_n \to u}
\int_0^t{\pd{}{t}\psi(s,u_n)\,ds} = \lim_{u_n \to u} \psi(t,u_n) -
u_n = \psi(t,u) - u\;.
\end{multline}
Since the left hand side of \eqref{Eq:R_ext_is_der} is
$t$-differentiable, also the right hand side is, and we obtain $R(u)
= \left.\pd{}{t}\psi(t,u)\right|_{t = 0}$ for all $u \in \cU$. A
similar calculation as above can be made upon replacing $R$ with
$F$, resulting in $F(u) = \left.\pd{}{t}\Phi(t,u)\right|_{t = 0}$
for all $u \in \cU$, and thus showing that the semi-homogeneous
affine process $\Xt$ is regular.
\end{proof}

\section{All affine processes are regular}\label{all}

In this final section we reduce the question of regularity of general
stochastically continuous affine processes to stochastically
continuous, \emph{semi-homogeneous} affine processes. Recall that for those
processes we have shown regularity in the preceding section.
The transformation of general affine processes to semi-homogeneous processes is
based on the method of the moving frame, which has been
successfully applied in the context of SPDEs several times; see for instance
\cite{filtaptei:08} and \cite{filtaptei:09}.

\begin{thm}
Every stochastically continuous affine process $ X $ is regular.
\end{thm}

\begin{proof}
By Theorem~\ref{Thm:Feller} $ X $ is a Feller process, and thus has a c\`adl\`ag version. Clearly, choosing a c\`adl\`ag version will not alter the functions $\Phi(t,u)$ and $\psi(t,u)$ defined by \eqref{Eq:Phi_definition}. Furthermore, we know
by Proposition \ref{Prop:J_exp} that
\begin{equation}
\psi_J(t,u) = \exp(t \beta) u_J
\end{equation}
for $(t,u) \in \cQ$ and a real $ n \times n $ matrix $ \beta $. We
define the $ d \times d $ matrix
\begin{equation}
K = \begin{pmatrix} \id_m & 0 \\ 0 & \beta
\end{pmatrix},
\end{equation}
and the transformation $\cT$
\begin{equation}\label{transformation}
Z_t = \cT[X]_t := X_t - K^\top\int_0^t{X_s\,ds},
\end{equation}
transforming the process $X$ path-by-path into a process $Z$. Note that the 
transformation is well-defined due to the c\`adl\`ag property of the
trajectories, and preserves the stochastic continuity of $ X $. Moreover, the transformation can be inverted by
\begin{equation}\label{inversion}
\cT^{-1}[Z]_t  = Z_t + K^\top \int_0^t {\exp \bigl((t-s)K^\top \bigr) Z_s\,ds},
\end{equation}
which is seen directly by inserting \eqref{transformation} and integrating by parts. 

We claim that the transformed process $Z = \cT[X]$ is a semi-homogeneous affine
process. For this purpose we calculate the conditional characteristic function:
Let $u \in i\RR^d$, and for each $N \in \NN$ and
$k \in \set{0,\dotsc,N}$, define $t_k = kt/N$, such that $t_0,
\dotsc, t_N$, is an equidistant partition of $[0,t]$ into
intervals of mesh $t/N$. By writing the time-integral as a limit
of Riemann sums, and using dominated convergence we have that
\begin{multline}\label{Eq:charf_Z}
\Excond{x}{\exp\left(\scal{u}{Z_{t+s}}\right)}{\cF_s} =
\exp\left(-\scal{u}{K^\top \int_0^s{X_r\,dr}}\right) \cdot \\
\lim_{N \to \infty}\Excond{x}{\exp\left(\scal{u}{X_{t+s}} -
\frac{t}{N}\scal{Ku}{\sum_{k=0}^{N-1}{X_{s + t_k}}}\right)}{\cF_s}\;.
\end{multline}
With the shorthands $h := t/N$ and $\Sigma_{n} := \sum_{k=0}^n{X_{s +
t_k}}$, and using the tower law as well as the affine property of
$\Xt$, the expectation on the right side can be written as
\begin{align*}
&\Excond{x}{\exp\Big(\scal{u}{X_{t+s}} - h \scal{Ku}{\Sigma_{N-1}}\Big)}{\cF_s}
= \\
&\Excond{x}{\vphantom{\Big(}\exp\Big(-h \scal{Ku}{\Sigma_{N-2}}\Big)
\cdot \Excond{x}{\exp\Big(\scal{(\id_d - hK)u}{X_{t+s}}\Big)}{\cF_{s
+
t_{N-1}}}}{\cF_s} = \\
&\Phi(h,(\id_d - hK)u) \cdot \\
& \qquad\cdot \Excond{x}{\exp\Big(\scal{\psi(h,(\id_d - hK)u)}{X_{s
+ t_{N-1}}} - h\scal{Ku}{\Sigma_{N-2}}\Big)}{\cF_s}\;.
\end{align*}
Applying the tower law $(N-1)$-times in the same way (conditioning
on \linebreak $\cF_{s + t_{N-1}}, \cF_{s + t_{N-2}},
\dotsc, \cF_{s + t_1}$, respectively) we arrive at the equation
\begin{align*}
& \Excond{x}{\exp\left(\scal{u}{X_{t+s}} -
\frac{t}{N}\scal{Ku}{\sum_{k=0}^{N-1}{X_{s + t_k}}}\right)}{\cF_s} \\
&= p(N-1;t,u) \exp\Big(\scal{X_s}{q(N-1;t,u)}\Big)\,,
\end{align*}
where the quantities $p(N-1;t,u)$ and $q(N-1;t,u)$ are defined through the
following recursion:
\begin{subequations}\label{Eq:recursion}
\begin{align}
&p(0;t,u) = 1, &\qquad& p(k+1;t,u) =
\Phi\Big(h,(\id_d - hK)q(k;t,u)\Big) \cdot p(k;t,u) \;,\\
&q(0;t,u) = u, &\qquad& q(k+1;t,u) = \psi\Big(h,(\id_d -
hK)q(k;t,u)\Big)\;,
\end{align}
\end{subequations}
Since the Riemannian sums in \eqref{Eq:charf_Z} converge point by point, we
conclude that the quantities $p(N-1;t,u) $ and $ q(N-1;t,u)$ converge to some
functions $p(t,u), q(t,u)$ as $N \to \infty$, and thus
\begin{equation}\label{Eq:Zt_indentity}
\Excond{x}{\exp\left(\scal{u}{Z_{t+s}}\right)}{\cF_s} = p(t,u)
\exp\left(\scal{q(t,u)}{X_s} -
\scal{u}{K^\top\int_0^s{X_r\,dr}}\right)
\end{equation}
for all $t,s \ge 0$ and $u \in i\RR^d$. Let now $q_J(t,u)$ denote
the $J$-components of $q(t,u)$. Based on the recursion
\eqref{Eq:recursion} and the fact that $\psi_J(t,u) = e^{\beta
t}u_J$ it holds that
\[q_J(t,u) = \lim_{N \to \infty} e^{t \beta}\Big(\id_n -
\frac{t\beta}{N}\Big)^{N-1} u_J = e^{t\beta}e^{-t\beta}u_J = u_J\;.\]
Thus we can rewrite \eqref{Eq:Zt_indentity} as
\[\Excond{x}{\exp\left(\scal{u}{Z_{t+s}}\right)}{\cF_s} =
p(t,u) \exp\left(\scal{q_I(t,u)}{Z^I_s} + \scal{u_J}{Z^J_s}\right)\]
which shows that $Z$ is indeed a stochastically continuous, semi-homogeneous
affine process.
By Theorem~\ref{Thm:condA_regular} such a process is regular. Hence, the functions $\wt{F}(u) = \left.\pd{}{t}p(t,u)\right|_{t = 0}$ and $\wt{R}(u) = \left.\pd{}{t}q(t,u)\right|_{t = 0}$ exist and satisfy the admissibility conditions in \citet[Def.~2.6]{Duffie2003}. By \citet[Thm.~2.7]{Duffie2003}, the functions $\wt{F}(u)$ and $\wt{R}(u) + K u$ are also admissible, and thus define a regular affine process $\wt{X}$. Using now the Feynman-Kac formula in \citet[Prop.~11.2]{Duffie2003}, it is seen that the transformation $\cT$ transforms the regular affine process $\wt{X}$ into the regular affine process $\cT[\wt{X}]$ characterized by $\wt{F}(u)$ and $\wt{R}(u)$ -- that is into a process equal in law to $Z$. We have shown that 
\[\cT[X] = Z \qquad \text{and} \qquad \cT[\wt{X}] = Z\,,\]
where equality is understood in law. Since the transformation $\cT$ can be inverted path-by-path, we conclude that $X = \wt{X}$ in law, and thus that $X$ is regular. 
\end{proof}

\begin{rem}
The intuition behind the `moving frame' transformation used above is the
following: given the process $ X $ we first construct a a time-dependent
coordinate transformation
$$
Y_t = \exp \bigl(- K^{\top} t \bigr) X_t,
$$
the `moving frame'. In the moving frame the process $X$ becomes time-dependent,
but can be re-scaled (in order to arrive at a time-homogeneous process) by the
stochastic integral
$$
d Z_t = \exp \bigl( K^{\top} t \bigr) d Y_t.
$$
The stochastic integral can be defined by integration by parts, i.e.,
$$
Z_t =  \exp \bigl( K^{\top} t \bigr) Y_t - K^\top \int_0^t \exp \bigl( K^{\top}
t \bigr) Y_r dr,
$$
which yields the transformation formula \eqref{transformation}. The method of
the moving frame is therefore an operation which allows to remove (or change)
the linear drift of an affine process.
\end{rem}

\begin{rem}
Now that we have shown that every stochastically continuous affine
process $X$ is regular, all the results of \citet{Duffie2003}
on regular affine processes apply to $X$. It follows in particular
that the set $\cQ$, introduced in \eqref{Eq:cO_def} is actually
equal to $\cU$, and hence simply connected. Thus, the
logarithm $\phi(t,u) = \log \Phi(t,u)$ is uniquely defined by
choosing the main branch of the complex logarithm, and we can write
\[\EE^x\left[e^{\scal{X_t}{u}}\right] =
\exp\left(\phi(t,u) + \scal{x}{\psi(t,u)}\right)\;,
\]
for all $(t,u) \in \cU$, as in \citet{Duffie2003}.
\end{rem}


\end{document}